\documentclass{amsart}

\usepackage{amsmath,amscd,amssymb,amsthm}
\usepackage{latexsym,enumerate,graphicx,geometry,hyperref,rotating,braket}
\usepackage[all]{xy}

\newcommand{\thefilename}{rcf}
\newtheorem{thm_}{Theorem}[section]
\newtheorem{lemma_}[thm_]{Lemma}
\newtheorem{eg_}[thm_]{Example}
\newtheorem{prop_}[thm_]{Proposition}
\newtheorem{def_}[thm_]{Definition}
\newtheorem{rk_}[thm_]{Remark}
\newtheorem{cor_}[thm_]{Corollary}
\newcommand{\thm}[1]{\begin{thm_}#1\end{thm_}}
\newcommand{\lemm}[1]{\begin{lemma_}#1\end{lemma_}}
\newcommand{\eg}[1]{\begin{eg_}#1\end{eg_}}
\newcommand{\prop}[1]{\begin{prop_}#1\end{prop_}}
\newcommand{\defi}[1]{\begin{def_}#1\end{def_}}
\newcommand{\rk}[1]{\begin{rk_}#1\end{rk_}}
\newcommand{\cor}[1]{\begin{cor_}#1\end{cor_}}
\newcommand{\pf}[1]{\begin{proof}#1\end{proof}}
\DeclareMathOperator{\Gal}{Gal}
\DeclareMathOperator{\Hom}{Hom}
\newcommand{\fracl}[3]{\genfrac{(}{)}{}{}{#1}{#2}_{#3}}
\newcommand{\fracn}[2]{\genfrac{(}{)}{}{}{#1}{#2}}

\newcommand{\card}[1]{\#({#1})}

\newcommand{\ZZ}{\mathbb Z}
\newcommand{\QQ}{\mathbb Q}

\renewcommand{\a}{\mathfrak a}
\renewcommand{\b}{\mathfrak b}
\renewcommand{\d}{\mathfrak d}
\newcommand{\f}{\mathfrak f}
\newcommand{\m}{\mathfrak m}
\renewcommand{\o}{\mathfrak o}
\newcommand{\p}{\mathfrak p}
\renewcommand{\P}{\mathfrak P}
\newcommand{\q}{\mathfrak q}
\renewcommand{\t}{\times}

\newcommand{\lr}{\longrightarrow}
\newcommand{\Llr}{\Longleftrightarrow}
\newcommand{\hr}{\hookrightarrow}
\newcommand{\eq}[1]{\begin{equation}#1\end{equation}}
\newcommand{\eqn}[1]{\begin{equation*}#1\end{equation*}}

\newcommand{\aln}[1]{\begin{align*}#1\end{align*}}

\newcommand{\enmt}[1]{\begin{enumerate}#1\end{enumerate}}
\renewcommand{\it}{\item}

\newcommand{\itm}[1]{\it[\upshape{(#1)}]}
\newcommand{\newnoindbf}[1]{\vspace{2mm}\noindent\textbf#1}

\begin{document}
\title[On Orders in Number Fields]{On Orders in Number Fields: Picard Groups, Ring Class Fields and Applications}
\author[C. Lv]{Chang Lv}
\address{Key Laboratory of Mathematics Mechanization\\
NCMIS, Academy of Mathematics and Systems Science\\
Chinese Academy of Sciences, Beijing 100190, P.R. China}
\email{lvchang@amss.ac.cn}
\thanks{The work of this paper was supported by the NNSF of China (Grants Nos. 61121062), 973 Project (2011CB302401) and the National Center for Mathematics and Interdisciplinary Sciences, CAS}
\author[Y. Deng]{Yingpu Deng}
\email{dengyp@amss.ac.cn}
\date{August 8, 2014}
\begin{abstract}
In this article, we focus on orders in arbitrary number fields, consider their Picard groups and finally obtain ring class fields corresponding to them. The Galois group of the ring class field is isomorphic to the Picard group. As an application, we give criteria of the solvability of the diophantine equation $p=x^2+ny^2$ over a class of imaginary quadratic fields where $p$ is a prime element.
\end{abstract}
\maketitle

\section{Introduction}\label{sec.intro}
Orders in number fields are widely used in many problems. A typical instance can be found in the book \cite{cox} by David A. Cox, where the author was dealing with the solvability of the diophantine equation $p=x^2+ny^2$ over $\ZZ$ for an odd prime $p$ and a positive integer $n$. Using the order $\ZZ[\sqrt{-n}]$ in imaginary quadratic field $\QQ(\sqrt{-n})$ along with the corresponding ring class field, Cox gave a beautiful criterion for the equation before, available for arbitrary $n>0$.

There are many other studies of ring class fields, say \cite{invitation, construct}. However, it seems that the result that there's an isomorphism between the Picard group of an order and the Galois group of it's corresponding ring class field only restricts to orders in imaginary quadratic fields, more generally, CM-fields (a totally imaginary quadratic extension of a totally real field). Here we consider orders in arbitrary number fields and their ring class fields, generalizing methods in \cite{cox}. By construction isomorphisms, we identify the Picard group as an generalized ideal class group, and therefore using the class field theory in the ideal-theoretic version we obtain a class field whose Galois group is isomorphic to the Picard group. Definitions and basic facts about orders are given in Section \ref{sec_facts}, and isomorphisms between groups are studied in Section \ref{sec_iso} subsequently. Then we can define the ring class fields in Section \ref{sec_rcf}. The last Section is dedicated to an application, considering the solvability of the diophantine equation $p=x^2+ny^2$ over a class of imaginary quadratic fields where $p$ is a prime element.

\section{Definitions and Basic Facts}\label{sec_facts}
We adopt the standard notations without additional explanation. If $R$ is a ring, we denote as $R^\t$ the unite group of $R$. Let $K$ be a number field, $\o_K$ the ring of integers of $K$, $\o$ an order in $K$. Denote the group of fractional ideals of $K$ by $J_K$, the principal fractional ideals by $P_K$ and the ideal class group by $Cl_K$. For the invertible and principal fractional ideals of $\o$ we use $J(\o)$ and $P(\o)$, respectively. We now give some definitions and notations about a modulus. One can see, say, \cite[\S IV.1]{janusz} for more details.
\defi{
A modulus in $K$ is a formal product $\m=\prod_\p\p^{n_\p}$ over all primes $\p$, finite or infinite, of $K$, where all integers $n_\p\ge0$ are nonzero only at finite many primes, zero at all complex primes and less than one at all real primes. We may write $\m=\m_0\m_\infty$ where $m_0$ is an $\o_K$-ideal and $m_\infty$ is a product of distinct real primes of $K$. If all $n_\p=0$ we set $m=1$.
}
Let
\aln{
K_\m&:=\set{\alpha/\beta|\alpha,\beta\in\o_K, \alpha\o_K,\beta\o_K\text{ relative prime to }\m_0},\\
K_{\m,1}&:=\set{\alpha\in K_\m|\alpha\equiv^*1\pmod\m}
}
where the notation $\alpha\equiv^*\beta\pmod\m$ for $\alpha,\beta\in K^\t$ means $v_p(\alpha\beta^{-1}-1)\ge n_\p$ (i.e. $\alpha\beta^{-1}\in1+\p^{n_\p}(\o_K)_\p$) for every $\p^{n_\p}\|\m_0$ and $\sigma(\alpha)>0$ for every infinite prime $\sigma$ dividing $\m_{\infty}$. Here $v_p$ is the normalized $\p$-adic exponential valuation and $(\o_K)_\p$ is the localization of $\o_K$ at $\p$. Let $J_K^\m$ be the subgroup of $J_K$ generated by all prime ideals not dividing $\m_0$, and $P_{K,1}^\m:=\set{\alpha\o_K|\alpha\in K_{\m,1}}$. We have basically the inclusions of groups $K_{\m,1}\subseteq K_\m\subseteq K^\t$ and $P_{K,1}^\m\subseteq J_K^\m$.

Now we recall (see \cite[\S I.12]{neukirch_alnt})
\defi{
$\notag$
\enmt{
\itm{a} Let $K/\QQ$ be an algebraic number field of degree $n$. An order of $K$ is a subring $\o$ of $\o_K$ which possesses an integral basis of length $n$. Clearly $\o_K$ is an order containing any other ones so it is called the maximal order.
\itm{b} An fractional ideal $\a$ of $\o$ is defined to be a finitely generated $\o$-submodule of $K$. $\a$ is called invertible if there exists a fractional ideal $\b$ such that $\a\b=\o$
\itm{c} Denote the group of invertible ideals of $\o$ by $J(\o)$. It contains the group
\eqn{
P(\o):=\set{\alpha\o|\alpha\in K^\t}
}
of principal ones. The group $Pic(\o):=J(\o)/P(\o)$ is called the Picard group of the order $\o$. In the case $\o=\o_K$, $Pic(\o_K)=Cl_K$ the class group of $K$.
\itm{d} The normalization of $\o$ is defined as the integral closure of $\o$ in $K$. In our case it is $\o_K$, and hence it is finitely generated over $\o$. We call the biggest ideal $\f$ of $\o_K$ contained in $\o$, in other words,
\eq{\label{eq_conductor}
\f=\set{a\in\o_K|a\o_K\subseteq\o},
}
the conductor of $\o$. Clearly $\f$ is also an $\o$-ideal. Since $o_K$ is f.g. over $\o$ we have $\f\neq0$.
\itm{e} A prime $\p\neq0$ of $\o$ is called regular if the localization $\o_\p$ is integrally closed. In fact the latter implying $\o_\p$ is a discrete valuation ring.
}
}
Let's list some basic facts about orders. For the proofs one can see \cite[\S I.12]{neukirch_alnt}. 
\prop{\label{prop_facts}
Let $\o$ be an order in $K$. Then we have
\enmt{
\itm{a} $\o$ is an one-dimensional noetherian integral domain with fractional field $K$.
\itm{b} A fractional ideal $\a$ of $\o$ is invertible iff $\a\o_\p$ is a fractional principal ideal of $\o_\p$ for every prime ideal $\p\neq0$.
\itm{c} Given a nonzero prime ideal $\p$ of $\o$, we have that $\p\nmid\f$ iff $\p$ is regular. If this is the case, then $\p\o_K$ is a prime ideal of $\o_K$. Here for a nonzero prime ideal $\p$ and an ideal $\a$ of $\o$, we denote $\a\subseteq\p$ by $\p\mid\a$, called $\p$ divides $\a$. Since $\o$ is one-dimensional, $\p$ is maximal, so $\p\mid\a$ iff $\p+\a\neq\o$.
\itm{d} The correspondence $\a\mapsto(\a\o_\p)_\p$ yields an isomorphism
\eq{\label{eq_jodecomp}
J(\o)\cong\bigoplus_\p P(\o_\p).
}
\itm{e} The group $\o_K^\t/\o^\t$, $(\o_K/\f)^\t$ and $(\o/\f)^\t$ are finite and one has
\eq{\label{eq_picocount}
\#Pic(\o)=\frac{h_K}{\card{\o_K^\t/\o^\t}}\frac{\#(\o_K/\f)^\t}{\#(\o/\f)^\t},
}
where $h_K=\#Cl_K$ is the class number of $K$.
}
}

\section{The Connection Between the Picard Group and the Generalized Class Group}\label{sec_iso}
In this section, we will construct isomorphisms in order to connect the Picard group to the generalized class group.
\lemm{\label{lem_invp}
A prime $\p$ of $\o$ is regular if and only if it is invertible.
}
\pf{
This is \cite[Exercise 5, \S I.12]{neukirch_alnt}.
}
\prop{\label{prop_ufjof}
Let $\a$ be an ideal of $\o$ relative prime to the conductor $\f$, i.e. $\a+\f=\o$. Then $\a$ is invertible and is uniquely factored into regular prime deals.
}
\pf{
Since $\o$ is noetherian so the set $S(\a):=\set{\p\text{ prime ideal of }\o|\a\subseteq\p}$ is finite (see, e.g. \cite[(3.4), \S I.3]{neukirch_alnt}). It follows that $\a\o_\q=\o_\q$ for $\q\not\in S(\a)$. But for $\p\in S(\a)$, $\p+\f=\o$ so $\f\not\subseteq\p$ i.e. $\p\nmid\f$. Then $\o_\p$ is discrete valuation ring by Proposition \ref{prop_facts} (c). So $\a\o_\p$ is a power of the maximal ideal $\p\o_\p$ of $\o_\p$, say, $e_\p$. Therefore by Proposition \ref{prop_facts} (b) and Lemma \ref{lem_invp} $\a$ and $\prod_{\p\in S(\a)}\p^{e_\p}$ are all in $J(\o)$ and clearly have the same image under the isomorphism \eqref{eq_jodecomp}. Hence $\a=\prod_{\p\in S(\a)}\p^{e_\p}$ and these $\p$'s are all regular. The uniqueness also follows from \eqref{eq_jodecomp}.
}
\lemm{\label{lem_bjptopok}
The extension and contraction of ideals give a bijection between $S(\o,\f)$ the set of prime ideals of $\o$ not dividing $\f$ and $S_K^\f$ the set of prime ideals of $\o_K$ not dividing $\f$.
}
\pf{
If $\p$ is a prime ideal of $\o$ and $\p\nmid\f$, i.e. $\p+\f=\o$, then $\p\o_K+\f=\o_K$ since $\f$ is ideal in both $\o$ and $\o_K$. Moreover, $\p\o_K$ is a prime ideal of $\o_K$ by Proposition \ref{prop_facts} (c). Conversely, if $\tilde\p$ is a prime ideal of $\o_K$ such that $\tilde\p\nmid\f$, i.e. $\f\not\subseteq\tilde\p$. Then $\p=\tilde\p\cap\o$ is a prime ideal of $\o$ and $\f\not\subseteq\p$, meaning $\p\nmid\f$. Now for all $\p\in S(\o,\f)$, $\p\subseteq\p\o_K\cap\o$ and $\p$ is maximal, and therefore $\p\o_K\cap\o=\p$. And similarly $\tilde\p=(\tilde\p\cap\o)\o_K$ for all $\tilde\p\in S_K^\f$ because $(\tilde\p\cap\o)\o_K$ is maximal. This completes the proof.
}
Let $\check J(\o,\f):=\set{\a\subseteq\o|\a+\f=\o}$ be the monoid of ideals of $\o$ not dividing $\f$. By Proposition \ref{prop_ufjof}, it is contained in $J(\o)$ and is a free abelian monoid over $S(\o,\f)$. Let $J(\o,\f)$ be the subgroup of $J(\o)$ generated by $\check J(\o,\f)$. We also denote by $\check J_K^\f$ the monoid of ideals of $\o_K$ not dividing $\f$. It generates $J_K^\f$. Then we have
\prop{\label{prop_jofisoikf}
$J(\o,\f)$ is a free abelian group over $S(\o,\f)$ and we have
\eq{\label{eq_jofisoikf}
J(\o,\f)\cong J_K^\f
}
and the correspondences are $\a\mapsto\a\o_K$ for $\a\in J(\o,\f)$ and $\tilde\a/\tilde\b\mapsto(\tilde\a\cap\o)/(\tilde\b\cap\o)$ for $\tilde\a,\tilde\b\in\check J_K^\f$.
}
\pf{
Since $\check J(\o,\f)$ and $\check J_K^\f$ are free abelian monoids over $S(\o,\f)$ and $S_K^\f$ respectively, the bijection in Lemma \ref{lem_bjptopok} clearly induces an isomorphism from $\check J(\o,\f)$ to $\check J_K^\f$, sending $\a$ to $\a\o_K$. Moreover, one can check the inverse map is the contraction of ideals. Then the isomorphism is canonically extended to an isomorphism from $J(\o,\f)$ to $J_K^\f$, abelian group generated by $\check J(\o,\f)$ and $\check J_K^\f$ respectively. This is the isomorphism we want. The fact that $J(\o,\f)$ is a free abelian group over $S(\o,\f)$ follows from it's definition and Proposition \ref{prop_ufjof}. And the correspondences are clear. 
}
Before we can proceed we shall note that in fact $P_{K,1}^\f$ has another expression:
\lemm{\label{lem_pk1f}
Let $\check P_{K,1}^\f:=\set{\alpha\o_K|\alpha\in\o_K, \alpha\equiv1\pmod{\f}}$. Then $P_{K,1}^\f$ is a subgroup of $J_K^\f$ generated by $\check P_{K,1}^\f$.
}
\pf{
Let $\tilde P_{K,1}^\f$ be the subgroup generated by $\check P_{K,1}^\f$. Let $\alpha/\beta\in K_\f$ with $\alpha,\beta\in\o$ relative prime to $\f$ and $\alpha/\beta\equiv^*1\pmod\f$. By $\alpha\o_K+\f=\o_K$ we can find $\gamma\in\o_K$ such that $\alpha\gamma\equiv1\pmod\f$, implying that $\alpha\gamma\equiv^*1\pmod\f$. But $\alpha/\beta=(\alpha\gamma)/(\beta\gamma)\equiv^*1\pmod\f$. It follows that $\beta\gamma\equiv^*1\pmod\f$, so $\beta\gamma\equiv1\pmod\f$ (see the properties of $\equiv^*$ in \cite[\S IV.1]{janusz}). This proves $P_{K,1}^\f\subseteq\tilde P_{K,1}^\f$. The other side is clear.
}
\rk{
In a similar manner we also know that
\eq{\label{eq_pk1f}
K_{\f,1}=\set{\alpha/\beta|\alpha,\beta\in\o_K,\text{ relative prime to }\f\text{ and }\alpha\equiv\beta\pmod\f}
}
}
For the following definition one can see \cite[\S V.6]{janusz}.
\defi{
A congruence subgroup of $K$ $($defined mod $\m$$)$ is a subgroup $H^\m$ of $J_K^\m$ such that $P_{K,1}^\m\subseteq H^\m\subseteq J_K^\m$. A basic fact is that the so called ray class group mod $\m$, $Cl_K^\m:=J_K^\m/P_{K,1}^\m$, is finite and so is $J_K^\m/H^\m$, called the generalized ideal class group.
}
Let $\check P(\o,\f):=\set{\alpha\o|\alpha\in\o,\alpha\o+\f=\o}\subseteq J(\o,\f)$, $P(\o,\f)$ the subgroup of $J(\o,\f)$ generated by $\check P(\o,\f)$ and set $Pic(\o,\f):=J(\o,\f)/P(\o,\f)$. Similarly let $\check P_{K,\o}^\f:=\set{\alpha\o_K|\alpha\in\o\text{ and }\alpha\o+\f=\o}\subseteq J_K$, $P_{K,\o}^\f$ the subgroup of $J_K$ generated by $\check P_{K,\o}^\f$. Then we have
\thm{\label{thm_picofisoclkof}
The group $P_{K,\o}^\f$ is a congruence subgroup defined mod $\f$. Moreover, let $Cl_{K,\o}^\f:=J_K^\f/P_{K,\o}^\f$. Then the isomorphism \eqref{eq_jofisoikf} in Proposition \ref{prop_jofisoikf} induces an isomorphism
\eq{\label{eq_picofisoclkof}
Pic(\o,\f)\cong Cl_{K,\o}^\f.
}
}
\pf{
Let $\alpha\in\check P_{K,1}^\f$ so that $\alpha\in\o_K$ with $\alpha\equiv1\pmod{\f}$. Since $\f\subseteq\o$ we have $\alpha\in1+\f\subseteq\o$ and $\alpha\o+\f=\o$. This implies $\check P_{K,1}^\f\subseteq\check P_{K,\o}^\f$ and hence $P_{K,1}^\f\subseteq P_{K,\o}^\f$. The inclusion $P_{K,\o}^\f\subseteq J_K^\f$ is trivial. So $P_{K,\o}^\f$ is a congruence subgroup defined mod $\f$. Now under the isomorphism $J(\o,\f)\cong J_K^\f$, we need only to verify that the image of $P(\o,\f)$ is $P_{K,\o}^\f$. But this is obvious by their definitions.
}
Our next goal is to connect $Pic(\o,\f)$ to $Pic(\o)$. Anyhow we have
\lemm{\label{lem_picofinjpico}
The inclusion $J(\o,\f)\subseteq J(\o)$ induces an monomorphism $Pic(\o,\f)\hr Pic(\o)$.
}
\pf{
Our task is to verify $P(\o)\cap J(\o,\f)=P(\o,\f)$. The side that $P(\o)\cap J(\o,\f)\supseteq P(\o,\f)$ is clear. To prove the other side, let $\alpha\o\in P(\o)$ where $\alpha\in K$. Suppose $\alpha\o\in J(\o,\f)$ so we can write $\alpha\o=\a/\b$ where $\a,\b\in\check J(\o,\f)$. Remember that $Pic(\o)$ is finite by \eqref{eq_picocount}, and then there exists positive integer $n$ such that $\b^n=\beta\o$ for some $\beta\in\o$. It follows that $\alpha\o=\a\b^{n-1}/\b^n$, implying that $\a\b^{n-1}=\alpha\beta\o$ and $\alpha\beta\in\o$. Note that $\b+\f=\o$ so $\beta\o+\f=\o$ and similarly $\alpha\beta\o+\f=\o$. Hence $\alpha\o=\alpha\beta\o/\beta\o\in P(\o,\f)$ and the proof completes.
}
In fact we find that $Pic(\o,\f)\cong Pic(\o)$. We will prove it by counting both sides. Since $Pic(\o)$ is finite, so is $Pic(\o,\f)$. Let $P_K^\f:=\set{\alpha\o_K|\alpha\in K_\f}$ and then we have the following short exact sequence:
\eq{\label{eq_ses_pj}
1\lr P_K^\f/P_{K,\o}^\f\lr J_K^\f/P_{K,\o}^\f\lr J_K^\f/P_K^\f\lr 1
}
By \eqref{eq_picofisoclkof}, $J_K^\f/P_{K,\o}^\f=Cl_{K,\o}^\f\cong Pic(\o,\f)$ so all the three groups above are finite. We now need a simple result:
\lemm{
For any modulus $\m$, there is a natural isomorphism
\eqn{
J_K^\m/J_K^\m\cap P_K\cong Cl_K
}
}
\pf{
We omit it here and one can see \cite[Corollary 1.5, \S IV.1]{janusz}.
}
By a similar but more easy argument as in the proof of Lemma \ref{lem_picofinjpico} we obtain that $J_K^\f\cap P_K=P_K^\f$. Hence we have
\eq{\label{eq_picofhk}
\#Pic(\o,\f)=\card{P_K^\f/P_{K,\o}^\f}h_K
}
by \eqref{eq_ses_pj}. Now we want to count $P_K^\f/P_{K,\o}^\f$. If we define
\eqn{
K_{\f,\o}:=\set{\alpha/\beta|\alpha,\beta\in\o\text{ relative prime to }\f}
}
then we have the following commutative diagram with exact rows:
\eqn{\xymatrix{
&1\ar[r] &\o_K^\t\cap K_{\f,\o}\ar[r]\ar@{_{(}->}[d] &K_{\f,\o}\ar[r]\ar@{_{(}->}[d] &P_{K,\o}^\f\ar[r]\ar@{_{(}->}[d] &1\\
&1\ar[r] &\o_K^\t\ar[r] &K_\f\ar[r] &P_K^\f\ar[r] &1
}}
the kernels of the three vertical arrows are trivial and hence by the snake lemma we have the exact sequence of the three cokernels:
\eq{\label{eq_ses_okp}
1\lr \o_K^\t/\o_K^\t\cap K_{\f,\o}\lr K_\f/K_{\f,\o}\lr P_K^\f/P_{K,\o}^\f\lr 1
}
where we consider the first tow groups. For the first one, $\o_K^\t/\o_K^\t\cap K_{\f,\o}$, we have the short exact sequence
\eq{\label{eq_ses_o}
1\lr \o_K^\t\cap K_{\f,\o}/\o^\t\lr \o_K^\t/\o^\t\lr \o_K^\t/\o_K^\t\cap K_{\f,\o}\lr 1
}
where $\o_K^\t/\o^\t$ is finite by Proposition \ref{prop_facts} (e). For the second one, $K_\f/K_{\f,\o}$, we know the correspondence $\alpha/\beta\mapsto \bar\alpha/\bar\beta$ gives a canonical epimorphism $K_\f\lr (\o_K/\f)^\t$ whose kernel is $K_{\f,1}$ by \eqref{eq_pk1f}. The case is similar for $K_{\f,\o}\lr (\o/\f)^\t$. Therefore we have the following commutative diagram with exact rows:
\eqn{\xymatrix{
&1\ar[r] &K_{\f,1}\cap K_{\f,\o}\ar[r]\ar@{_{(}->}[d] &K_{\f,\o}\ar[r]\ar@{_{(}->}[d] &(\o/\f)^\t\ar[r]\ar@{_{(}->}[d] &1\\
&1\ar[r] &K_{\f,1}\ar[r] &K_\f\ar[r] &(\o_K/\f)^\t\ar[r] &1
}}
The snake lemma again yields a exact sequence
\eq{\label{eq_ses_k}
1\lr K_{\f,1}/K_{\f,1}\cap K_{\f,\o}\lr K_\f/K_{\f,\o}\lr (\o_K/\f)^\t/(\o/\f)^\t\lr 1
}
Now we have three short exact sequences \eqref{eq_ses_okp}, \eqref{eq_ses_o} and \eqref{eq_ses_k}. One can check all groups in them are finite and we can count $Pic(\o,\f)$ according these exact sequences and the relation \eqref{eq_picofhk}:
\aln{
\#Pic(\o,\f)&=h_K\card{P_K^\f/P_{K,\o}^\f}=h_K\frac{\card{K_\f/K_{\f,\o}}}{\card{\o_K^\t/\o_K^\t\cap K_{\f,\o}}}\\
&=h_K\frac{\card{K_{\f,1}/K_{\f,1}\cap K_{\f,\o}}\card{(\o_K/\f)^\t/(\o/\f)^\t}}{\card{\o_K^\t/\o^\t}/\card{\o_K^\t\cap K_{\f,\o}/\o^\t}}\\
&=\frac{h_K}{\card{\o_K^\t/\o^\t}}\frac{\#(\o_K/\f)^\t}{\#(\o/\f)^\t}\card{K_{\f,1}/K_{\f,1}\cap K_{\f,\o}}\card{\o_K^\t\cap K_{\f,\o}/\o^\t}.
}
Comparing to \eqref{eq_picocount}, we have $\#Pic(\o,\f)=\#Pic(\o)\card{K_{\f,1}/K_{\f,1}\cap K_{\f,\o}}\card{\o_K^\t\cap K_{\f,\o}/\o^\t}\ge\#Pic(\o)$. Recall we have the monomorphism $Pic(\o,\f)\hr Pic(\o)$ in Lemma \ref{lem_picofinjpico}. Hence the equality holds and we have
\thm{\label{thm_picofisopico}
There is a natural isomorphism
\eq{\label{eq_picofisopico}
Pic(\o,\f)\cong Pic(\o).
}
}
Also by our counting for groups we have
\prop{
We have $K_{\f,1}\subseteq K_{\f,\o}$ and $\o_K^\t\cap K_{\f,\o}=\o^\t$.
}

\section{The Ring Class Field Corresponding to the Order}\label{sec_rcf}
Let's recall a result of class field theory, stated in the ideal-theoretic version:
\thm{[The existence theorem]\label{thm_cft_existance}
Let $\m$ be a modulus of $K$, $H^\m$ a congruence subgroup defined mod $\m$. Then there exists a unique abelian extension $L$ of $K$ such that all primes of $K$ ramified in $L$, finite or infinite, divide $\m$, and if
\eqn{
\varphi_{L/K}^\m: J_K^\m\lr \Gal(L/K)
}
is the Artin map of $L/K$, then $H^\m=\ker\varphi_{L/K}^\m$. Moreover, if it is the case, the Artin reciprocity law holds for $(L,K,\m)$, i.e. $\ker\varphi_{L/K}^\m=P_{K,1}^\m N_{L/K}(J_L^\m)$. Since the Artin map is sujective, we also have
\eqn{
J_K^\m/H^\m=J_K^\m/P_{K,1}^\m N_{L/K}(J_L^\m)\cong \Gal(L/K).
}
}
\pf{
One can see \cite[Theorem 5.8, 9.9 and 11.11, Chapter V]{janusz}.
}
At this point we have all we need to obtain the ring class field. This is our main
\thm{[The ring class field]\label{thm_rcf}
Let $\o$ be an order of $K$ with conductor $\f$. Then there exists a unique abelian extension $H_\o$ of $K$, such that all primes of $K$ ramified in $H_\o$ $($must be finite$)$ divide $\f$, and the kernel of the Artin map $\varphi_{H_\o/K}^\f:J_K^\f\lr \Gal(H_\o/K)$ is $P_{K,\o}^\f$. At this time we have the induced isomorphism
\eqn{
Pic(\o)\cong \Gal(H_\o/K).
}
We call $H_\o$ the ring class field corresponding to the order $\o$.
}
\pf{
This is a direct corollary to Theorem \ref{thm_picofisoclkof}, \ref{thm_picofisopico} and \ref{thm_cft_existance}.
}
\rk{
It is obvious that in the case $\o=\o_K$, $H_\o$=$H_K$ the Hilbert class field of $K$.
}
\cor{\label{cor_rcf_split}
Let notations be as in Theorem \ref{thm_rcf}. Let $\tilde\p$ be a prime ideal of $\o_K$ that is relative prime to $\f$. Then $\tilde\p$ splits completely in $H_\o$ if and only if $\tilde\p=\alpha\o_K$ for some $\alpha\in\o$.
} 
\pf{
Let $\p=\tilde\p\cap\o$. Since $\tilde\p\in J_K^\f$ we have $\p\in J(\o,\f)$. We first claim that $\tilde\p\in P_{K,\o}^\f$ iff $\tilde\p=\alpha\o_K$ for some $\alpha\in\o$. The isomorphism $Pic(\o,\f)\cong Cl_{K,\o}^\f$ sends $\p$ to $\tilde\p$. Thus if $\tilde\p\in P_{K,\o}^\f$ then $\p\in P(\o,\f)\cap\o=\check P(\o,\f)$, so $\p=\alpha\o$ for some $\alpha\in\o$ and hence $\tilde\p=\p\o_K=\alpha\o_K$. The converse of the claim is trivial. Then the result follows from the properties that the Artin map sends a unramified prime ideal to identity of the Galois group iff the prime ideal splits completely.
}

\section{An Application to Diophantine Equations}
The main theorem of David A. Cox \cite{cox} is a beautiful criterion of the solvability of the diophantine equation $p=x^2+ny^2$. The detail is
\thm{\label{thm_cox_2h}
Let $n>0$ be an integer. Then there is a monic irreducible polynomial $f_n(x)\in\ZZ[x]$ of degree $h(-4n)$ such that if an odd prime $p$ divides neither $n$ nor the discriminant of $f_n(x)$, then $p=x^2+ny^2$ is solvable over $\ZZ$ if and only if $\fracl{-n}{p}{}=1$ and $f_n(x)=0$ is solvable over $\ZZ/p\ZZ$. Here $h(-4n)$ is the class number of primitive positive definite binary forms of discriminant $-4n$. Furthermore, $f_n(x)$ may be taken to be the minimal polynomial of a real algebraic integer $\alpha$ for which $L=K(\alpha)$ is the ring class field of the order $\ZZ[\sqrt{-n}]$ in the imaginary quadratic field $K=\QQ(\sqrt{-n})$.
}
There are generalizations by Dasheng Wei \cite{wei1, wei2} using results on the Manin obstruction, giving a criterion of the solvability of the $\alpha=x^2+y^2$, over the integer rings of a class of quadratic fields. That is for arbitrary $\alpha$ but the only case $n=1$. In this section, we also give some generalizations over a class of imaginary quadratic fields for various $n$ but still prime element $p$, as an application of the generalization of ring class fields obtained in the previous sections. The idea is also a generalization of that in \cite{cox}.

\newnoindbf{{Notations and definitions.}
We fix some notations that will be used in this section. From now  on let $F$ be a number field, $\o_F$ the ring of integers of $F$, $n$ an element in $\o_F$ such that $-n$ is not a square in $F$. Let $E=F(\sqrt{-n})$, and $\o_E$ the corresponding ring of integers. Note that $\Gal(E/F)=\ZZ/2\ZZ$, and by writing a bar we mean the action by the non-trivial element of $\Gal(E/F)$ sending $\sqrt{-n}$ to $-\sqrt{-n}$. So $\p\o_E=\P\bar\P$ for some $\P$ prime ideal of $E$ if a prime ideal $\p$ of $F$ splits completely in $E$. Let $\d(E/F)$ be the discriminant of $E/F$ and $\d(\alpha)$ be the discriminant of an element $\alpha\in E$.
Set $\o=\o_F+\o_F\sqrt{-n}$ and clearly this is an order of $E$. Still denote the ring class field corresponding to $\o$ by $H_\o$. We are going to see that this order has some special properties. 
\prop{\label{prop_conductor}
The conductor $\f$ of $\o$ contains $4n$.
}
\pf{
By the definition of the conductor \eqref{eq_conductor} we only need to show that $4n\o_E\in\o$. Let $\alpha=a+b\sqrt{-n}\in\o_E$ with $a,b\in F$. We may assume $\alpha\not\in F$, whence the minimal polynomial of $\alpha$, i.e. $x^2-2ax+a^2+nb^2$, is in $\o_F[x]$. Hence $2a\in\o_F$ and $(4nb)^2=4n(4(a^2+nb^2)-(2a)^2)\in\o_F$. This means $4nb\in F$ is integral over $\o_F$ and then $4nb\in\o_F$. It follows that $4n\alpha=(2n)(2a)+(4nb)\sqrt{-n}\in\o_F+\o_F\sqrt{-n}=\o$ and we have done.
} 
Note that the only non-trivial element of $\Gal(E/F)$ (denoted as a bar) sends $\o$ to itself isomorphically, so by the definition of $\f$ we have $\bar\f=\f$.
\prop{\label{prop_rcf_gal}
The ring class field $H_\o$ is Galois over $F$.
}
\pf{
Let $\sigma\in\Hom_F(H_\o, F^{\text{al}})$ be an $F$-embedding from $H_\o$ into the algebraic closure of $F$, it suffices to show that $\sigma H_\o=H_\o$ for any such $\sigma$. Clearly $\sigma|_E\in\Gal(E/F)=\ZZ/2\ZZ$ so $\sigma H_\o$ is also an abelian extension over $\sigma E=E$. We have seen that $\bar\o=\o$ and $\bar\f=\f$. Hence $\sigma\o=\o$ and $\sigma\f=\f$. It follows that $\sigma P_{E,\o}^\f=P_{E,\o}^\f$ and the Artin map of $\sigma H_\o/E$ is defined on $I_E^\f$. Recall that Theorem \ref{thm_rcf} asserts $H_\o$ is the only abelian extension of $E$ such that $\ker\varphi_{H_\o/E}^\f=P_{E,\o}^\f$. By the observation that $\varphi_{\sigma H_\o/E}^\f(\sigma\a)=\sigma\varphi_{H_\o/E}^\f(\a)\sigma^{-1}$ for $\a\in I_E^\f$, we know that $\ker\varphi_{\sigma H_\o/E}^\f=\sigma P_{E,\o}^\f=P_{E,\o}^\f$. Hence $\sigma H_\o$ is also an extension with this property. Therefore $\sigma H_\o=H_\o$ and the proof is complete.
}

\subsection{Criteria by Solutions of Polynomials}
Given a nonzero prime ideal $\p$ of $\o_F$, we are able to give a rough criterion of the solvability of the equation
\eq{\label{eq_ideal}
\p=(x^2+ny^2)\o_F
}
of ideals.
}
\prop{\label{prop_basic}
Let $\p$ be a nonzero prime ideal of $\o_F$. If $2n\not\in\p$ then there exist $x$ and $y$ in $\o_F$, such that $\p=(x^2+ny^2)\o_F$ if and only if $\p$ splits completely in $H_\o$, the ring class field corresponding to the order $\o=\o_F+\o_F\sqrt{-n}$.
}
\pf{
Since $2n\not\in\p$ and $\d(E/F)\mid\d(\sqrt{-n})=-4n$ we have $\p\nmid \d(E/F)$. Hence $\p$ is unramified in $E$ (see, e.g. \cite[Theorem 7.3]{janusz}). Proposition \ref{prop_conductor} implies $\p\o_E\nmid\f$ and then we shall prove the following equivalences:
\aln{
\p=(x^2+ny^2)\o_F&\Llr \p\o_E=\P\bar\P\text{ where }\P=\alpha\o_E\text{ prime to }\f\text{ for some }\alpha\in\o\\
&\Llr \p\o_E=\P\bar\P\text{ where }\P\subset\o_E\text{ is prime to }\f\text{ and splits completely in }H_\o\\
&\Llr \p\text{ splits completely in }H_\o.
}
and the theorem will follow. Now if $\p=(x^2+ny^2)\o_F=(x+\sqrt{-n}y)(x-\sqrt{-n}y)\o_F$, writing $\P=(x+\sqrt{-n}y)\o_E$ and $\alpha=x+\sqrt{-n}y\in\o$, we have $\p\o_E=\P\bar\P$ where $\P=\alpha\o_E$. Moreover, since $\bar\f=\f$ and $\P\not=\bar\P$ it follows that $\p\o_E\nmid\f$ if and only if $\P\nmid\f$. In fact, this is the prime decomposition of $\p$ in $\o_E$. Conversely, suppose that $\p\o_E=\P\bar\P$ with $\P=\alpha\o_E$ where $\alpha\in\o$. Hence $\P=(x+\sqrt{-n}y)\o_E$ for some $x$ and $y$ in $\o_F$. It follows that $\p\o_E=(x^2+ny^2)\o_E$, implying $\p=(x^2+ny^2)\o_F$ by unique factorization of ideals in $\o_E$. This proves the first equivalence. The second equivalence comes from Corollary \ref{cor_rcf_split}. And the converse of the last equivalence is trivial. To see the other side, we just use the fact that $\P\not=\bar\P$, $[E:F]=2$, and that $H_\o/F$ is Galois (Proposition \ref{prop_rcf_gal}). The proof is complete.
}
We now derive two criteria of the ideal equation \eqref{eq_ideal} involving solvabilities of polynomials over residue class fields. The cases become easier since the residue class fields are all finite fields.
\thm{\label{thm_2h}
Let $\p$ be a nonzero prime ideal of $\o_F$. Then there is a monic irreducible polynomial $g_n(x)\in\o_F[x]$ of degree $2h(\o)$, such that if $\p$ contains neither $2n$ nor the discriminant of $g_n(x)$, then $\p=(x^2+ny^2)\o_F$ for some $x$ and $y$ in $\o_F$ iff $g_n(x)=0$ is solvable over $\o_F/\p$. Here $h(\o)=\#Pic(\o)$.
}
\pf{
Let $H_\o=F(\beta)$ with $\beta\in\o_{H_\o}$ and $g_n(x)\in\o_F[x]$ be the minimal polynomial of $\beta$, which clearly has degree $2h(\o)$ by Theorem \ref{thm_rcf}. If $\p$ does not contain the discriminant of $g_n(x)$, which is $\d_{H_\o/F}(\beta)$, we could use Kummer's Lemma (see \cite[Proposition 25, \S I.8]{gtm110}, along with Proposition 16, Section III.3 in the same book) and remember that $H_\o/F$ is Galois. It follows that $\p$ splits completely in $H_\o$ iff $g_n(x)=0$ is solvable over $\o_F/\p$. Hence the theorem follows by Proposition \ref{prop_basic}.
}

\subsection{Applied to a Class of Imaginary Quadratic Fields}
In this subsection, we consider imaginary quadratic fields. Let $d>3$ and $n$ be two distinct positive square free integers. Let $F=\QQ(\sqrt{-d})$. So we have $E=F(\sqrt{-n})=\QQ(\sqrt{-d}, \sqrt{-n})$, a biquadratic field. We will make some further restrictions on $d$ and $n$ in order to make the ideal equation \eqref{eq_ideal} into the ordinary equation
\eq{
p=x^2+ny^2,
}
mainly by considering the units. We give the following criterion.
\thm{\label{thm_quadr_2h}
Let $F=\QQ(\sqrt{-d})$ where $d>3$ and $n$ are two distinct positive square free integers. Assume that
\eq{\label{eq_unit}
-1=\alpha^2+n\beta^2\text{ for some }\alpha,\beta\in \o_F.
}
Then there is a monic irreducible polynomial $g_n(x)\in\o_F[x]$ of degree $2h(\o)$, such that if a prime element $p\in\o_F$ divides neither $2n$ nor the discriminant of $g_n(x)$, then $p=x^2+ny^2$ is solvable over $\o_F$ iff $g_n(x)=0$ is solvable over $\o_F/p\o_F$.
}
\pf{
We use Theorem \ref{thm_2h} to obtain that $p\o_F=(x^2+ny^2)\o_F$ if and only if $g_n(x)=0$ is solvable over $\o_F/p\o_F$. The former is clearly implied by $p=x^2+ny^2$. For the other side, since $d>3$ the units of $\o_F$ are $\{\pm1\}$. It follows that if $p\o_F=(x^2+ny^2)\o_F$ then $p=\pm(x^2+ny^2)$. If the plus sign holds, we have done; otherwise by the assumption \eqref{eq_unit}, then we have $p=(\alpha^2+n\beta^2)(x^2+ny^2)=(\alpha x-n\beta y)^2+n(\alpha y+\beta x)^2$. So we complete the proof.
}
In fact the assumption \eqref{eq_unit} is true for a lot of cases. For one of these classes of imaginary quadratic fields we have the
\cor{\label{cor_quadr_speci_2h}
Let $F=\QQ(\sqrt{-d})$ where $d\ge5$ and $n\ge1$ are tow distinct rational primes of the form $4k+1$. If $\fracn{d}{n}=-1$ or $\fracl{d}{n}{4}=-1$ then there is a monic irreducible polynomial $g_n(x)\in\o_F[x]$ of degree $2h(\o)$, such that if a prime element $p\in\o_F$ divides neither $2n$ nor the discriminant of $g_n(x)$, then $p=x^2+ny^2$ is solvable over $\o_F$ iff $g_n(x)=0$ is solvable over $\o_F/p\o_F$.
}
\pf{
By Theorem \ref{thm_quadr_2h} it suffices to prove that the assumption \eqref{eq_unit} holds. Since $d$ and $n$ are of the form $4k+1$, if $\fracn{d}{n}=-1$, or $\fracn{d}{n}=1$ and $\fracl{d}{n}{4}=\fracl{n}{d}{4}=-1$, then $x^2-dny^2=-1$ is solvable over $\ZZ$ (\cite[p. 228]{dirichlet}); otherwise $\fracn{d}{n}=1$, $\fracl{d}{n}{4}=-1$ and $\fracl{n}{d}{4}=1$, so $x^2-dny^2=d$ is solvable over $\ZZ$ (\cite[Corollary 1.4]{wei_diophantine}). Both cases indicate that there exists $\alpha,\beta\in\o_F$ such that $-1=\alpha^2+n\beta^2$. The proof is done.
}
It should be noted that there is no general but only \emph{ad hoc} method to determine the defining polynomial of the generalized ring class field (i.e. the minimal polynomial of an integral primitive element of $H_\o/E$)  For this reason, we could not turn Corollary \ref{cor_quadr_speci_2h} into explicit criteria.

\subsection{Finding Hilbert Class Fields}\label{sec.hilb}
From now no we consider the case where
\eq{\label{eq_o_eq_ok}
\o=\o_E\quad\text{ i.e. }\quad\o_E=\o_F+\o_F\sqrt{-n}
}
so that we are able to give explicit criteria for some cases. Note that at this time $H_\o$ is the Hilbert class field and we denote it as $H_E$.
In this case a lower degree polynomial can be used instead of $g_n(x)$, just as in \cite{cox}. This can be done when there is a defining polynomial of $H_E/E$ which happens to be in $\o_F[x]$.
\thm{\label{thm_hilb_h}
Let $\p$ be a nonzero prime ideal of $\o_F$ and assume \eqref{eq_o_eq_ok} holds. Suppose further that there is a field $K$ such that $H_E=KE$ and $F=K\cap E$. Then there is a monic irreducible polynomial $f_n(x)\in\o_F[x]$ of degree $h(E)$, the class number of $E$, such that if $\p$ contains neither $2n$ nor the discriminant of $f_n(x)$ then
\aln{
&\p=(x^2+ny^2)\o_F\text{ for some }x,y\in\o_F\\
\Llr\quad &-n\text{ is a square mod }\p\text{ and }f_n(x)=0\text{ is solvable over }\o_F/\p.
}
}
\pf{
The proof is similar to Theorem \ref{thm_2h}. Since we already have by the proof of Proposition \ref{prop_basic} that
\aln{
&\p=(x^2+ny^2)\o_F\text{ for some }x,y\in\o_F\\
\Llr\quad &\p\o_E=\P\bar\P\text{ , and }\P\text{ splits completely in }H_E,
}
by Kummer's Lemma, noting that $\p\nmid4n$, $\p\o_E=\P\bar\P$ iff $x^2+n=0$ is solvable in $\o_F/\p$, which is to say $-n$ is a square mod $\p$. Now consider the following diagram of field extensions
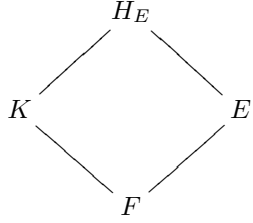
\eqn{\xymatrix{
&	&H_E \ar@{-}[dl] \ar@{-}[dr]\\
&K \ar@{-}[dr]	&	&E \ar@{-}[dl]\\
&	&F
}}
Suppose $K=F(\alpha)$ with $\alpha\in\o_K$. Since $H_E$ is the composition $KE$ we have $H_E=E(\alpha)$. Let $f_n(x)\in\o_F[x]$ be the minimal polynomial of $\alpha$ over $F$. By $K\cap E=F$ we know that $[K:F]=[H_E:E]=h(E)$. Hence we obtain that $f_n(x)$ is also the minimal polynomial of $\alpha$ over $E$, of degree $h(E)$. If $\p$ does not divide the discriminant of $f_n(x)$, using Kummer's Lemma again, then $\P$ splits completely in $H_E$ iff $f_n(x)=0$ is solvable in $\o_E/\P$. Remember we already know $\p=\P\bar\P$ splits completely in $E$, so $\o_E/\P\cong\o_F/\p$, and hence $f_n(x)=0$ is solvable in $\o_E/\P$ iff it is solvable in $\o_F/\p$. This completes the proof.
}
There are results on ensuring the existence of $K$ for some cases. One is given by
\thm{[Wyman]\label{thm_K}
Let $E/F$ be a cyclic extension and suppose that $h(F)=1$. Then there exists $K$ such that $H_E=KE$ and $F=K\cap E$.
}
\pf{
This is taken from \cite[Proposition 3]{splt.hilb}.
}
Hence we have the
\cor{\label{cor_hilb_K_h}
Let $\p$ be a nonzero prime ideal of $\o_F$ and assume \eqref{eq_o_eq_ok} holds. Suppose further that $h(F)=1$. Then there is a monic irreducible polynomial $f_n(x)\in\o_F[x]$ of degree $h(E)$, the class number of $E$, such that if $\p$ contains neither $2n$ nor the discriminant of $f_n(x)$ then
\aln{
&\p=(x^2+ny^2)\o_F\text{ for some }x,y\in\o_F\\
\Llr\quad &-n\text{ is a square mod }\p\text{ and }f_n(x)=0\text{ is solvable over }\o_F/\p.
}
}
\pf{
We have $h(F)=1$, and noting that $E/F$ is Galois, thus Theorem \ref{thm_K} ensures the existence of $K$, and Theorem \ref{thm_hilb_h} applies.
}
Next we will give another ensuring condition, leading not only the existence but also what $K$ can be taken to be.

If we use the above criteria in Theorem \ref{thm_hilb_h} for computation, we must attempt to find the corresponding Hilbert class fields (or the field $K$). One can use complex multiplication (see, e.g. \cite{cox}) for those of quadratic imaginary fields, Stark's method for those of totally real number fields, and Kummer theory for those of a general number field (see \cite{gtm193}). The main focus of this section is on solvability over imaginary quadratic fields, which in turn, needs Hilbert class fields of biquadratic fields (see the beginning of this section). So we will give some constructions of Hilbert class fields in special cases, in fact the constructions of defining polynomials. We first need a
\prop{[E. Friedman]
Let $E/F$ be number fields. Let $Cl(E)$ and $Cl(F)$ denote the class group of $E$ and $F$, respectively. Consider the exact sequence
\eqn{\xymatrix{
&0 \ar[r]	&Cl_N(E/F) \ar[r]	&Cl(E) \ar[r]^{N_{E/F}}	&Cl(F) \ar[r]	&Cl_{N, E}(F) \ar[r]	&0
}}
where $Cl_N(E/F)$ and $Cl_{N, E}(F)$ denote the kernel and cokernel of the norm map $N_{E/F}$. Let $H_E$ and $H_F$ denote the corresponding Hilbert class fields. Then we have the following assertions.
\enmt{
\it[\upshape{(a)}] $Cl_N(E/F)\cong \Gal(H_E/EH_F)$. In particular, $N_{E/F}$ is injective iff $EH_F=H_E$.
\it[\upshape{(b)}] $N_{E/F}$ is injective iff $[E\cap H_F:F]=\card{Cl(F)}/\card{Cl(E)}$.
\it[\upshape{(c)}] $Cl_{N, E}(F)\cong \Gal(E\cap H_F/F)$. In particular, $N_{E/F}$ is surjective iff $E\cap H_F=F$.
}
}
\pf{
This is taken from \cite[Exercise 3, \S7.6]{gtm193}.
}
By this proposition,
\eq{\label{eq_inj}
[E\cap H_F:F]=\card{Cl(F)}/\card{Cl(E)}
}
if and only if $EH_F=H_E$. This helps us to find $H_E$ in some cases, just by finding $H_F$ first and then verifying \eqref{eq_inj}. Another application is the more special case that
\eq{\label{eq_inj1}
[E\cap H_F:F]=\card{Cl(F)}/\card{Cl(E)}=1,
}
which is equivalent to that the norm map is an isomorphism or to that $EH_F=H_E$ and $E\cap H_F=F$. The later means that we can take $K=H_F$ in Theorem \ref{thm_hilb_h}, even if Theorem \ref{thm_K} dos not apply. We summarize it as the following
\thm{\label{thm_hilb}
Assume that the assumptions \eqref{eq_o_eq_ok} and \eqref{eq_inj1} hold. Let $f(x)\in\o_F$ $($of degree $h(E)=h(F)$$)$ be a defining polynomial of $H_F/F$ $($which is independent on $n$$)$. If $p\in\o_F$ a prime element divides neither $2n$ nor the discriminant of $f(x)$, then $p\o_F=(x^2+ny^2)\o_F$ for some $x$ and $y$ in $\o_F$ iff $-n$ is a square mod $p$.
}
\pf{
Since the prime ideal $\p=p\o_F$ is principal, it splits completely in $H_F$. The later is to say that $f(x)=0$ is solvable over $\o_F/p\o_F$ since $\p$ does not divide the discriminant of $f(x)$. Applying Theorem \ref{thm_hilb_h}, taking $K=H_F$ and $f_n(x)=f(x)$ there, we have
\aln{
&p\o_F=(x^2+ny^2)\o_F\text{ for some }x,y\in\o_F\\
\Llr\quad &-n\text{ is a square mod }p\text{ and }f(x)=0\text{ is solvable over }\o_F/p\o_F.
}
However it already holds that $f(x)=0$ is solvable over $\o_F/p\o_F$. The proof completes.
}
\rk{
By considering modulo $p$, one knows clearly that $p\o_F=(x^2+ny^2)\o_F$ implies $-n$ is a square mod $p$.
}
We see from this Theorem that if \eqref{eq_inj1} holds, the criterion becomes very simple, just whether $-n$ is a square mod $p$. In particular, if
\eq{\label{eq_injo}
\card{Cl(F)}=\card{Cl(E)}\text{ is odd,}
}
then \eqref{eq_inj1} holds because the extension degrees of $E$ and $H_F$ over $F$ are relative prime.

Now we are going to deal with quadratic fields. Still let $d>3$ and $n$ be two distinct positive square free integers and let $F=\QQ(\sqrt{-d})$.
We will make some further restrictions on $d$ and $n$ in order to satisfy the assumption \eqref{eq_o_eq_ok}.
\cor{\label{cor_hilb_quadr_speci}
Let $d>3$ and $n$ be two positive square free integers. Suppose they are relative prime and $d\equiv3\pmod4$ and $n\equiv1\text{ or }2\pmod4$. Assume the diophantine equation $du^2-nv^2=1$ is solvable over $\ZZ$. Suppose further that \eqref{eq_inj1} holds. In particular, suppose \eqref{eq_injo} holds. Let $f(x)\in\o_F$ be a defining polynomial of $H_F/F$. If $p$ divides neither $2n$ nor the discriminant of $f(x)$, then $p=x^2+ny^2$ is solvable over $\o_F$ iff $-n$ is a square mod $p$.
}
\pf{
Since we assume that $d\equiv3\pmod4$ and $n\equiv1\text{ or }2\pmod4$, \cite[Exercise 42(c), Chapter 2]{marcus} tells us that
\eqn{
\{1, \frac{1+\sqrt{-d}}{2}, \sqrt{-n}, \frac{1+\sqrt{-d}}{2}\sqrt{-n}\}
}
is an integral basis for $\o_E$. It follows that $\o_E=\o_F+\o_F\sqrt{-n}$, i.e. the assumption \eqref{eq_o_eq_ok} holds. Noting the assumption \eqref{eq_inj1} we can use Theorem \ref{thm_hilb} to obtain that $p\o_F=(x^2+ny^2)\o_F$ for some $x$ and $y$ in $\o_F$ iff $-n$ is a square mod $p$. The condition $du^2-nv^2=1$ implies that \eqref{eq_unit} holds. Like the proof of Theorem \ref{thm_quadr_2h} we know $p\o_F=(x^2+ny^2)\o_F$ is equivalent to $p=x^2+ny^2$.
}
Let us exhibit an example for the case in Corollary \ref{cor_hilb_quadr_speci}. The computation of class numbers is done by the package \emph{PARI} \cite{pari} mainly developed by Henri Cohen.
\eg{
$d=59$ and $n=2$, $F=\QQ(\sqrt{-59})$, if a prime element $p$ of $\o_F=\ZZ[\frac{1+\sqrt{-59}}{2}]$ dos not divide $118$ then $p=x^2+2y^2$ is solvable over $\o_F$ iff $-2$ is a square mod $p$.
}
\pf{
$F=\QQ(\sqrt{-59})$ and $E=F(\sqrt{-2})$. $h(F)=h(E)=3$ is odd so \eqref{eq_injo} holds. Hence Corollary \ref{cor_hilb_quadr_speci} applies. We need to find a defining polynomial of $H_F/F$. We know by Theorem \ref{thm_K} (applied to $F/\QQ$) that there is a integer coefficient one. We refer to the table for Hilbert class fields of several imaginary quadratic fields \cite[Appendix 12.1.2]{gtm193}, from which we know $x^3+2x-1$ is what we want. In fact we may also use complex multiplication to find this. Since the discriminant of $x^3+2x-1$ is $-59$, the result follows.
}
We continue this example in details. For instance, consider the rational prime $17$. It is easy to see that $17$ splits completely in $F$ since $\fracn{-59}{17}=1$. In fact $17=p\o_F\overline{p\o_F}$ where $p=\frac{3+\sqrt{-59}}{2}$ a prime element in $\o_F$. Clearly $p\nmid118$. Since here $n\in\ZZ$ and $\o_F/p\o_F\cong\ZZ/17\ZZ$, we find that the criterion for $p=x^2+2y^2$ being solvable over $\o_F$ in the example becomes that $\fracn{-2}{17}=1$. This clearly holds. Hence it is true that $p=x^2+2y^2$ is solvable over $\o_F$. In fact we have
\eqn{
\frac{3+\sqrt{-59}}{2}=(\frac{5779+1115\sqrt{-59}}{2})^2+2(-3028+266\sqrt{-59})^2.
}

\newnoindbf{{Acknowledgment}
The authors would like to thank Pace Nielsen, Yupeng Jiang, Dandan Huang and Caihua Luo for many helpful discussions and comments. Special thanks to Xintong Zhao for her revision on the writing.
}

\bibliography{\thefilename}

\providecommand{\bysame}{\leavevmode\hbox to3em{\hrulefill}\thinspace}
\providecommand{\MR}{\relax\ifhmode\unskip\space\fi MR }
\providecommand{\MRhref}[2]{%
  \href{http://www.ams.org/mathscinet-getitem?mr=#1}{#2}
}
\providecommand{\href}[2]{#2}
\begin{thebibliography}{10}

\bibitem{pari}
\emph{{PARI/GP}}, \url{http://pari.math.u-bordeaux.fr/}.

\bibitem{gtm193}
Henri Cohen, \emph{Advanced topics in computational number theory}, Graduate
  Texts in Mathematics, vol. 193, Springer, 2000.

\bibitem{construct}
Harvey Cohn, \emph{Introduction to the construction of class fields}, Courier
  Dover Publications, 1985.

\bibitem{invitation}
Harvey Cohn and Taussky Olga, \emph{A classical invitation to algebraic numbers
  and class fields}, Springeri, New York, 1978.

\bibitem{splt.hilb}
Gary Cornell and Michael Rosen, \emph{A note on the splitting of the {Hilbert}
  class field}, Journal of Number Theory \textbf{28} (1988), no.~2, 152--158.

\bibitem{cox}
David~A. Cox, \emph{Primes of the form {$x^2+ ny^2$}: Fermat, class field
  theory, and complex multiplication}, John Wiley \& Sons, 1989.

\bibitem{dirichlet}
G.L. Dirichlet, \emph{Einige neue s\"atze \"uber unbestimmte gleichungen},
  "Werke" \textbf{I} (1920), 221--236.

\bibitem{janusz}
Gerald~J. Janusz, \emph{Algebraic number fields}, 2 ed., Graduate Studies in
  Mathematics, vol.~7, American Mathematical Soc., 1996.

\bibitem{gtm110}
Serge Lang, \emph{Algebraic number theory}, Graduate Texts in Mathematics, vol.
  110, Springer, New York, 1970.

\bibitem{marcus}
Daniel~A Marcus, \emph{Number fields}, Springer, 1977.

\bibitem{neukirch_alnt}
J\"urgen Neukirch, \emph{Algebraic number theory}, Springer, 1999.

\bibitem{wei1}
Dasheng Wei, \emph{On the sum of two integral squares in quadratic fields
  {$\mathbb Q(\sqrt{\pm p})$}}, Acta Arith. \textbf{147} (2011), no.~3,
  253--260.

\bibitem{wei_diophantine}
\bysame, \emph{On the {Diophantine} equation {$x^2-Dy^2=n$}}, Science China
  Mathematics \textbf{1-12} (2013).

\bibitem{wei2}
\bysame, \emph{On the sum of two integral squares in the imaginary quadratic
  field {$\mathbb Q(\sqrt{-2p})$}}, Science China Mathematics \textbf{57}
  (2014), no.~1, 49--60.

\end{thebibliography}
\bibliographystyle{amsplain} 
\end{document}